\documentclass[uslettersize, 10pt, conference]{ieeeconf}      

\IEEEoverridecommandlockouts                              
\usepackage{graphics} 
\usepackage{epsfig} 
\usepackage{mathptmx} 
\usepackage{amsmath} 
\usepackage{amssymb}  
\usepackage{amscd}
\usepackage{ifthen}
\usepackage{cite}

\def\scr#1{{\cal #1}}
\newcommand{\R}{\mathbb{R}}
\def\eq#1{\begin{equation}#1\end{equation}}
\def\rep#1{(\ref{#1})}

\newtheorem{theorem}{Theorem}

\newtheorem{remark}{Remark}

\title{
A Stackelberg Game for Multi-Period Demand Response Management \\in the
Smart Grid
\thanks{K. Alshehri, J. Liu, X. Chen, and T. Ba\c{s}ar are with the Coordinated Science Laboratory, University of Illinois at Urbana-Champaign (\texttt{\{kalsheh2, jiliu, xdchen, basar1\}@illinois.edu}).
}
}

\author{Khaled Alshehri, Ji Liu, Xudong Chen, and Tamer Ba\c sar}

\begin{document}

\maketitle
\thispagestyle{empty}
\pagestyle{empty}

\begin{abstract}
This paper studies a multi-period demand response management problem in the smart grid where multiple utility companies compete among themselves. The user-utility interactions are modeled by a noncooperative game of a Stackelberg type where the interactions among the utility companies are captured through a Nash equilibrium.
It is shown that this game has a unique Stackelberg equilibrium
at which the utility companies
set prices to maximize their revenues (within a Nash game) while the users respond accordingly to maximize their utilities
subject to their budget constraints.
Closed-form expressions are provided for the corresponding strategies of the users and the utility companies.
It is shown that the multi-period scheme, compared with the single-period case,
provides more incentives for the users to participate in the game.
A necessary and sufficient condition on the minimum budget needed for a user to participate
is provided.
\end{abstract}

\section{Introduction}

%
%

Demand Side Management (DSM) is an essential component of the smart grid as it captures important aspects of the interactions between utility companies (UCs) and consumers, including residential, commercial, industrial users and vehicles \cite{survey}. These aspects can be technical, such as using advanced metering infrastructure to improve the reliability and efficiency of the grid \cite{survey}, or social through agreements between consumers and energy providers \cite{DOECOM}. The categorization of DSM includes energy efficiency, time-of-use pricing, and demand response \cite{DSM}. Energy efficiency is about the long-term solutions to improve the overall efficiency; one example is replacing the air conditioning system of a building with a more efficient one. Through time-of-use pricing, companies can agree with the users that different prices are to be charged at different periods; these prices are agreed on beforehand by both parties. Demand Response Management (DRM)  is the response of consumers' demands to price signals from the UCs \cite{DR}. DRM allows companies to manage the consumers' demands, either directly or indirectly, through incentive-based programs \cite{review}.

Game theory is a powerful mathematical tool that can lead to effective multi-person decision making \cite{basar}. Due to the nature of the smart grid which features different  entities with conflicting objectives, applying game-theoretic methods  can improve their reliability and efficiency \cite{survey,ss,trading,amir,walidPHEV,repeated,sabita,fourstage,hazem,walidPHEV2,twolevel}. For a comprehensive survey of game-theoretic methods for communications, DSM, and microgrid systems in the smart grid, we refer to \cite{survey}. 

There are several works where game-theoretic methods have been applied to DSM/DRM \cite{amir,walidPHEV,repeated,sabita,fourstage,hazem,walidPHEV2,twolevel}.
An autonomous DSM through scheduling of appliances has been implemented within a noncooperative game framework in \cite{amir}. The participants in such a game are energy users who are connected to the same UC and the outcome of the game is the power consumption schedule of appliances that minimizes the overall energy cost. Plug-in hybrid electric vehicles can sometimes have a surplus of energy, and it can be more beneficial to sell energy quickly instead of holding onto it. This motivates the work of \cite{walidPHEV} in which an energy trading noncooperative game between vehicle groups and distribution grids has been studied. In \cite{repeated}, joint consumers discomfort and billing costs minimization within a repeated game framework has been shown to lead to an optimal DSM mechanism.

  In the smart grid, there are numerous interactions between UCs and the end-users because of the two-way communication infrastructure \cite{DOECOM}. In order to take this fact into account, utilizing a multi-level framework such as Stackelberg games would be very useful \cite{sabita,fourstage,hazem,walidPHEV2}. A multi-level game-theoretic framework has been developed for demand response management, where consumers choose their optimal demands in response to prices announced by  different UCs \cite{sabita}. This Stackelberg game was shown to have a unique Stackelberg equilibrium at which UCs maximize their revenues and end-users maximize their payoff functions. In this framework, UCs were the leaders of the game and users were the followers.  This framework, though it effectively captures user-utility interactions, is limited to the single-period case. In the smart grid, users are expected to be able to schedule their energy consumption, store or sell surplus energy, based on their self-interests. With energy scheduling and storage, users might have a flexibility about when to receive a certain amount of energy, particularly for shiftable appliances \cite{amir}. Some energy consumption can be scheduled and some cannot. For example, users might be flexible about doing the laundry, but not much flexibility is there for refrigerators.

   It is worth mentioning that multi-level games for DRM have been studied in a limited context in the literature. For example, a four-stage Stackelberg game has been studied where three stages are at the leader-level (the utility retailer), and the fourth stage is at the consumer level \cite{fourstage}. Retailers choose the amount of energy to procure, and the sources to produce it, in addition to deciding on the price. Consumers respond to these prices through demand selection. This game is also a single-period game and it does not take into account the competition between the UCs, but, it incorporates other aspects of the decision making at the company-level. Furthermore, a multi-period Stackelberg game with a single UC has been formulated with energy storage devices at the consumers-level \cite{hazem}.  Simulation results show that a Peak-to-Average Ratio minimization problem leads to similar decisions made within the framework of the game. Additionally, a noncooperative Stackelberg game between plug-in electric vehicle groups such as parking lots and the smart grid was formulated and a socially optimal equilibrium has been achieved \cite{walidPHEV2}. A two-level game (a noncooperative game between multiple UCs and an evolutionary game for the users at the lower level) has been proposed in\cite{twolevel}, but this game is also limited to the single-period case.

  There is a need to develop an analytical multi-period, multi-utility, and a multi-level framework. By introducing multi-period inter-temporal constraints, we could study a generalization of \cite{sabita} to the multi-period case. Our work differs from the work in \cite{sabita} at both the user-side and the company-side. At the user-side, we have an additional minimum energy constraint that needs to be satisfied across all periods, while at the company-side, we provide an alternative computationally cheap closed-form solution for the prices. Having such a multi-period framework can make it possible to accommodate numerous extensions, such as, energy scheduling and storage, and Peak-to-Average ratio minimization. To recap, the multi-period scheme has been considered before in the literature, but most of the studies are either user-oriented or UC-oriented. Our goal here is to develop a reliable multi-period framework where numerous multi-level interactions are also taken into account.

Accordingly, we formulate in this paper a Stackelberg game for multi-period-multi-company demand response management. We derive solutions in closed-form and find precise expressions for the maximizing demands at the users' level, and the revenue-maximizing prices for the UCs. We also prove the existence and uniqueness of the Stackelberg equilibrium. This work captures the competition between UCs and the budget limitations at the consumer-level.

The remainder of this paper is organized as follows.
The problem is formulated in Section \ref{formulation}. The Stackelberg game along with its solutions in closed form are provided in Section \ref{game}. In Section \ref{simulation}, some numerical results are provided. We conclude the paper in Section \ref{conclusion} with some remarks on future directions.

\section{Problem Formulation} \label{formulation}

The system consists of $K\geq 1$ UCs and $N \geq 1$ end-users.
We consider a finite time horizon including $T\geq1$ periods. Let $\scr{K}=\{1,2,\dots,K\}$ be the set of UCs, $\scr{N}=\{1,2,\dots,N\}$ be the set of end-users, and $\scr{T}=\{1,2,\dots,T\}$ be the finite set of time slots.
The unit of time can be hour, day, week, or month.
Mathematically, it does not matter which unit the set $\scr{T}$ represents.

We model the interaction between UCs and end-users as a Stackelberg game. Thus, the DRM problem
is formulated as two optimization problems, one is for the  followers of the game and the other is for the leaders of the game. The followers of the game are the consumers, who can be residential, commercial, or industrial, and the leaders of the game are the UCs.
On the followers (users) side, each user has a certain amount of energy demand to meet and we assume that the users have flexibility to schedule their consumption of energy. It is also assumed that users have limited budgets. The goal of each user is to maximize her utility without exceeding the budget while meeting the minimum energy demand. On the leaders (UCs) side, the goal of each UC is to maximize its revenue while its price setting is influenced by both the competition against other UCs and the behavior of the users. In the game considered here, the UCs announce their prices for each period to the users, and then the users respond accordingly by scheduling their demands. We discuss the underlying models for both sides separately below.
\subsection{User-Side}
Because of energy scheduling and storage, users may have some flexibility about when to receive a certain amount of energy. In our analysis, we are concerned about the total amount of shiftable energy. For non-shiftable energy, one can add some period-specific constraints.  Additionally, each user has a budget constraint that she cannot exceed for the entire time horizon.

  The demand of user $n\in \scr{N}$ from UC $k\in \scr{K}$ at time $t\in \scr{T}$ is denoted by $d_{n,k}(t)$, and  $p_k(t)$ denotes the price of energy announced by UC $k$ at time $t$. Each end-user $n$ has a budget of $B_n$, and  a minimum amount of energy to be met, denoted by $E^{{\rm min}}_n$. The value of $E_n^{{\rm min}}$ can be thought of as the total energy needed by user $n$ during the considered total period (for example, one day or one week). User payoff functions increase by receiving more energy.
The utility of user $n$ is defined as \begin {equation}{ U_{{\rm user},n}=\gamma_n\sum_{k\in \scr{K}}\sum_{t\in \scr{T}}\ln(\zeta_n+d_{n,k}(t))} \label{user}\end{equation}
where $\gamma_n$ and $\zeta_n$ (typically, $\zeta_n=1$) are adjustment parameters. The logarithmic function is well known to provide a good demand response \cite{sabita,basarDR}.
Users aim to maximize their payoffs as much as possible while meeting this threshold of minimum amount of energy and not exceeding a certain budget. To be more precise,
given $B_n \geq 0$, $E^{{\rm min}}_n \geq 0$, and $p_k(t)>0$, the user-side optimization problem is defined as follows:
\begin{eqnarray}
\underset{\mathbf{d_{n,k}}}{\hbox{maximize}} && U_{{\rm user},n} \nonumber \\
\hbox{subject to} && \sum_{k\in \scr{K}}\sum_{t\in \scr{T}}p_k(t)d_{n,k}(t)\leq B_n \label{cc1}\\
&& \sum_{k\in \scr{K}}\sum_{t\in \scr{T}}d_{n,k}(t)\geq \,E^{{\rm min}}_n \label{cc}\\
&&d_{n,k}(t)\geq 0,\;\;   \forall k \in \scr{K},\;\;  \forall t \in \scr{T}
\end{eqnarray}
Note that there is no game played among the users. Each user responds to the price signals using only her local information. These price signals depend on all the demands selected by the users and hence users indirectly affect each other's decisions.
\subsection{Company-Side}

Given the prices of other UCs and the power availability of UC $k$ at period $t$, denoted by ${\bf{p_{-k}}}$ and $G_k(t)$ respectively, the total revenue for UC $k$ is given by
\begin{equation}{U_{{\rm gen},k}(p_k,{\bf{p_{-k}}})=\sum_{t\in \scr{T}}p_k(t)\sum_{n\in \scr{N}}d_{n,k}(t)}\label{UC}\end{equation}
The utility-side optimization problem is then described as follows:
\begin{eqnarray}
\underset{\mathbf{p_k}}{\hbox{maximize}} && U_{{\rm gen},k} (p_k,{\bf{p_{-k}}})
\nonumber \\
\hbox{subject to} && \sum_{n\in \scr{N}} d_{n,k}(t) \leq G_k(t),\;\; \forall \; t \in \scr{T} \label{prob2} \\
&& p_k(t)> 0, \;\; \forall \; t \in \scr{T}
\end{eqnarray}

  The goal of each UC is to maximize its revenue and hence maximize its profit.  Also, companies want to meet the consumers' demands without exceeding the maximum power availability (overloading the system can cause contingencies). Additionally, because of the market competition, the prices announced by other companies also affect the determination of the price at company $k$. So, company $k$'s price selection is actually a response to what other UCs in the market have announced; this response is also constrained by the availability of power. Thus, there is a Nash game between utility companies.

\section{The Stackelberg Game} \label{game}

In this section, we first solve the above optimization problems in closed-form. Then, we study the existence and uniqueness of Nash equilibrium at the UC-level and the Stackelberg equilibrium for the entire two-level game.

\subsection{Followers-Side Analysis}
Note that the user-side utility function  is convex and the constraints are linear.
We start by relaxing the minimum energy constraint (\ref{cc}) and then find the necessary budget that makes the maximizing demands feasible.
  For each user $n\in \scr{N}$, the associated Lagrange function is given as follows:
\begin{eqnarray*}
L_{{\rm user},n} &=& \gamma_n\sum_{k\in \scr{K}}\sum_{t\in \scr{T}}\ln(\zeta_n+d_{n,k}(t))
 \\
&&- \lambda_{n,1}\left(\sum_{k\in \scr{K}}\sum_{t\in \scr{T}}p(t)_kd_{n,k}(t)-B_n\right) \\
&&+\sum_{k\in\scr{K}} \sum_{t\in\scr{T}} \lambda_{n,2}(k,t)d_{n,k}(t)
\end{eqnarray*}
  where $\lambda_{n,i}$'s are the Lagrange multipliers. For optimality, by Krush-Kuhn-Tucker necessary conditions \cite{NL},
\begin{eqnarray}
&&\frac{\partial L_{{\rm user},n}}{\partial d_{n,k}(t)} =0,\;\; \forall \;t\in \scr{T}, \;\; \forall \;k\in \scr{K}\\
&&\lambda_{n,1}\left(\sum_{k\in \scr{K}}\sum_{t\in \scr{T}}p_k(t)d_{n,k}(t)-B_n \right) = 0 \label{c1} \\
&&\lambda_{n,2}(k,t) d_{n,k}(t) = 0, \;\;\forall \;t\in \scr{T}, \;\;\forall \;k\in \scr{K}  \label{c2}\\
&&\lambda_{n,1},\lambda_{n,2}(k,t) \geq 0, \;\; \forall \;t\in \scr{T}, \;\;\forall \; k\in \scr{K}
\end{eqnarray}
The above conditions are also sufficient because of the convexity of the optimization problem \cite{boyd}.
In the sequel,
we derive the solution in closed form by discussing two cases.

\subsubsection{Case 1:  $d_{n,k}(t)>0, k\in \scr{K}, t\in \scr{T}$}

In this case, since $\lambda_{n,1}>0$  (constraint (\ref{cc1}) is active \cite{NL}) and $\lambda_{n,2}(k,t)=0$
for all $t\in \scr{T}$ and $k\in \scr{K}$, there holds
\eq{\frac{\partial L_{{\rm user},n}}{\partial d_{n,k}(t)} = \frac{\gamma_n}{\zeta_n+d_{n,k}(t)}-\lambda_{n,1}p_k(t)=0,
\;\forall \;t\in \scr{T}, \;k\in \scr{K}
\label{x}}
From (\ref{x}) and (\ref{c1}), it is straightforward to verify that
\begin{equation}
d_{n,k}(t)= \frac{B_n+\sum_{k\in \scr{K}}\sum_{t\in \scr{T}}p_k(t)\zeta_n}{KTp_k(t)}-\zeta_n,  \; \forall \;t\in \scr{T}, \; k\in \scr{K}  \label{xx}
\end{equation}
which is a generalization of the single-period case considered in \cite{sabita}.

\subsubsection{Case 2: at least one $d_{n,k}(t)$ is zero}

We show that the expression (\ref{xx}) also holds for $d_{n,k}(t)\geq0$, $k\in \scr{K}$, $t\in \scr{T}$. Without loss of generality, suppose that $d_{n,1}(1)=0$ and $d_{n,e}(f)>0$ for $ e\in \scr{K}$ and $f\in \scr{T}$, except when $(e,f)=(1,1)$. Following a similar analysis as in the previous case,
$$d_{n,e}(f)= \frac{B_n+\sum_{e\in \scr{K}}\sum_{f\in \scr{T}}p_e(f)\zeta_n}{(KT-1)p_e(f)}-\zeta_n$$
Note that since  $d_{n,1}(1)=0$, there holds
$$B_n=\zeta_n\left(KTp_1(1)-\sum_{k\in \scr{K}}\sum_{t\in \scr{T}}p_k(t)\right)$$
Thus, it is straightforward to verify that
$$d_{n,e}(f)=
\frac{B_n+\sum_{k\in \scr{K}}\sum_{t\in \scr{T}}p_k(t)\zeta_n}{KTp_e(f)}-\zeta_n$$
which matches the expression (\ref{xx}). By the budget constraint (i.e., $\sum_{k\in \scr{K}}\sum_{t\in \scr{T}}p_k(t)d_{n,k}(t)\leq B_n $), we can see that when $d_{n,k}(t)=0$ for all $k\in \scr{K}$ and $t\in \scr{T}$, user $n$ has  a zero budget (i.e., $B_n=0$). The announced prices cannot be infinite because of the nature of the Stackelberg equilibrium, as we discuss later. The following theorem states the necessary and sufficient condition for $B_n$ so that the above $d_{n,k}(t)$'s are guaranteed to be feasible.

\begin{theorem}
For each user $n \in \scr{N}$, the demand $d_{n,k}(t)$ given by (\ref{xx}) is feasible if and only if $$B_n \geq max \{f_{n,1},f_{n,2}\}$$ where
\begin{eqnarray*}
f_{n,1} &=& \zeta_n\left(KTp_k(t)-\sum_{k\in \scr{K}}\sum_{t\in \scr{T}}p_k(t)\right)\\
f_{n,2} &=& \frac{E_n^{{\rm min}}+\zeta_nKT}{\sum_{k\in \scr{K}}\sum_{t\in \scr{T}}\frac{1}{KTp_k(t)}}-\zeta_n \sum_{k\in \scr{K}}\sum_{t\in \scr{T}}p_k(t)
\end{eqnarray*}
\end{theorem}

\vspace{.2in}
\begin{proof}
Suppose first that $B_n \geq f_{n,1}$. Then,
$$B_n \geq \zeta_n(KTp_k(t)-\sum_{k\in \scr{K}}\sum_{t\in \scr{T}}p_k(t)) \,\,\,\, \forall k \in \scr{K},\,\, t \in \scr{T}$$
which leads to
$$\frac{B_n+\zeta_n\sum_{k\in \scr{K}}\sum_{t\in \scr{T}}p_k(t)}{KTp_k(t)} -\zeta_n\geq 0 \,\,\,\, \forall k \in \scr{K},\,\, t \in \scr{T}$$
Thus, the inequality $B_n \geq f_{n,1}$ implies that the non-negativity condition is satisfied.
Next suppose that $B_n \geq f_{n,2}$. Then,
$$B_n \geq \frac{E_n^{min}+\zeta_nKT}{\sum_{k\in \scr{K}}\sum_{t\in \scr{T}}\frac{1}{KTp_k(t)}}-\zeta_n \sum_{k\in \scr{K}}\sum_{t\in \scr{T}}p_k(t)$$
Thus,
$$B_n+\zeta_n \sum_{k\in \scr{K}}\sum_{t\in \scr{T}}p_k(t) \geq \frac{E_n^{min}+\zeta_nKT}{\sum_{k\in \scr{K}}\sum_{t\in \scr{T}}\frac{1}{KTp_k(t)}}$$
from which we have
$$\sum_{k\in \scr{K}}\sum_{t\in \scr{T}}\frac{B_n+\zeta_n \sum_{k\in \scr{K}}\sum_{t\in \scr{T}}p_k(t)}{KTp_k(t)}- \sum_{k\in \scr{K}}\sum_{t\in \scr{T}}\zeta_n \geq E_n^{min}$$
This leads to
$$\sum_{k\in \scr{K}}\sum_{t\in \scr{T}}d_{n,k}(t)\geq \,E^{min}_n$$
With this, the inequality $B_n \geq f_{n,2}$ implies that the minimum energy need is satisfied.
Combining both conditions, we conclude that the condition
$$B_n \geq \max \{f_{n,1},f_{n,2}\}$$
guarantees that the maximizing demand in (\ref{xx}) is feasible.

Now suppose that
$$\frac{B_n+\zeta_n\sum_{k\in \scr{K}}\sum_{t\in \scr{T}}p_k(t)}{KTp_k(t)} -\zeta_n \geq 0 \,\,\,\, \forall k \in \scr{K},\,\, t \in \scr{T}$$
It follows that
$$B_n \geq \zeta_n(KTp_k(t)-\sum_{k\in \scr{K}}\sum_{t\in \scr{T}}p_k(t)) \,\,\,\, \forall k \in \scr{K},\,\, t \in \scr{T}$$
which implies that $B_n \geq f_{n,1}$.
Finally, suppose that
\begin{eqnarray*}
\sum_{k\in \scr{K}}\sum_{t\in \scr{T}}d_{n,k}(t) &=& \sum_{k\in \scr{K}}\sum_{t\in \scr{T}}(\frac{B_n+\zeta_n \sum_{k\in \scr{K}}\sum_{t\in \scr{T}}p_k(t)}{KTp_k(t)} - \zeta_n) \\
&\geq& E^{min}_n
\end{eqnarray*}
Then,
$$\sum_{k\in \scr{K}}\sum_{t\in \scr{T}}\frac{B_n+\zeta_n \sum_{k\in \scr{K}}\sum_{t\in \scr{T}}p_k(t)}{KTp_k(t)} \geq E^{min}_n+KT\zeta_n$$
By re-arranging, we have
$$B_n \geq \frac{E_n^{min}+\zeta_nKT}{\sum_{k\in \scr{K}}\sum_{t\in \scr{T}}\frac{1}{KTp_k(t)}}-\zeta_n \sum_{k\in \scr{K}}\sum_{t\in \scr{T}}p_k(t) = f_{n,2}$$
from which it follows that  $B_n \geq \max \{f_{n,1},f_{n,2}\}$ when the demand is feasible.
\end{proof}

\subsection{Leaders-Side Analysis}

Given the prices set by the other companies and subject to the power availability constraint (\ref{prob2}), each UC (leader) aims to determine its most profitable prices. At the leaders level, there is a noncooperative game in which each UC chooses its optimal prices in response to the prices set by the other UCs. The revenue function of company $k$ is an increasing function of the consumers' demands $\sum_{n\in \scr{N}} d_{n,k}(t)$. Thus, the optimality is reached when the equality in (\ref{prob2}) holds. We apply the solutions derived in the users-side analysis (which was a function of the prices) here and show that in the case when the equality in (6) is satisfied, the corresponding prices constitutes the best response of UC $k$ subject to the prices set by the other UCs.

With the equality in (\ref{prob2}) and relation (\ref{xx}), there holds
\begin{eqnarray*}
&&\frac{\sum_{n\in \scr{N}}B_n+\sum_{n\in \scr{N}}\zeta_n\sum_{k\in \scr{K}}\sum_{t\in \scr{T}}p_k(t)}{KTp_k(t)}\\
&&=\sum_{n\in \scr{N}}\zeta_n + G_k(t),\;\; \forall \;t \in \scr{T}\\
\end{eqnarray*}
Let $Z=\sum_{n\in \scr{N}}\zeta_n$ and $B=\sum_{n\in \scr{N}}B_n$. Then, for each company $k \in \scr{K}$,
$$ B+Z\sum_{k\in \scr{K}}\sum_{t\in \scr{T}}p_k(t) = KTp_k(t)(G_k(t)+Z),\;\; \forall \;t \in \scr{T}$$
Note that the double summation includes $p_k(t)$ and all the other prices.
Thus, \begin{equation}\begin{split}B+Z\sum_{e\in \scr{K}}\sum_{h\in \scr{T}}p_e(h)= KTp_k(t)(G_k(t)+Z)-p_k(t)Z, \\ \,\,\,\, \forall \;t \in \scr{T}, \;\forall \;k \in \scr{K}, \; (e,h)\neq (k,t) \label{AP}\end{split}\end{equation}
Note that the equations in (\ref{AP}) can be combined into a linear equation
$$AP= Y$$
where $A$ is a $KT\times KT$ matrix whose diagonal entries are $KT(G_k(t)+Z)-Z$, $k\in\scr{K}$, $t\in\scr{T}$,
and off-diagonal entries all equal $-Z$,
$P$ is a vector in $\R^{KT}$ stacking $p_k(t)$, $k\in\scr{K}$, $t\in\scr{T}$,
and $Y$ a vector in $\R^{KT}$ whose entries all equal $B$.

The following theorem shows that matrix $A$ is invertible and the revenue-maximizing prices are positive and unique.
\begin{theorem}
The following statements are true.
\begin{enumerate}
\item The matrix $A$ is nonsingular and the prices announced by company $k\in \scr{K}$ at time $t\in \scr{T}$ are uniquely given by
    \begin{equation} p_k(t)=\frac{B}{G_k(t)+Z}\left(\frac{1}{KT-\sum_{k\in \scr{K}}\sum_{t\in \scr{T}}\frac{Z}{G_k(t)+Z}}\right)\label{p}\end{equation}
\item The above prices  are always positive.
\item The price given by \rep{p} constitutes the best response of company $k$ to the prices set by other companies.
\end{enumerate}
\end{theorem}

\vspace{.1in}
\begin{proof}
\begin{enumerate}
\item The matrix $A$ can be represented as
\begin{footnotesize}
\begin{equation*}\begin{split} A=\begin{pmatrix}
KT(G_1(1)+Z)& 0&\dots& 0\\
0 & KT(G_1(2)+Z)& \dots& 0\\
\vdots & \ddots\\
0 &\dots&0 & KT(G_K(T)+Z)
\end{pmatrix}\\+\begin{pmatrix}
-Z\\
-Z\\
\vdots \\
-Z
\end{pmatrix}\begin{pmatrix}
1 \dots 1
\end{pmatrix}:=\hat{A}+uv^T\end{split}\end{equation*}
\end{footnotesize}

Note that $\hat{A}$ is invertible (diagonal matrix with nonzero diagonal elements). By Sherman-Morrison Formula \cite{SM},  if $\hat{A}$ is invertible and $1+v^T\hat{A}^{-1}u\neq0$, then $$(\hat{A}+uv^T)^{-1}=\hat{A}^{-1}-\frac{\hat{A}^{-1}uv^T\hat{A}^{-1}}{1+v^T\hat{A}^{-1}u}$$
Note that
\begin{footnotesize}
\begin{align}\nonumber
&1+v^T\hat{A}^{-1}u=1-\cr
&\begin{pmatrix}
1 \dots 1
\end{pmatrix}\begin{pmatrix}
\frac{1}{KT(G_1(1)+Z)}& 0&\dots& 0\\
0 & \frac{1}{KT(G_1(2)+Z)}& \dots& 0\\
\vdots & \ddots\\
0 &\dots&0 & \frac{1}{KT(G_K(T)+Z)}
\end{pmatrix}  \begin{pmatrix}
Z\\
Z\\
\vdots \\
Z
\end{pmatrix} \cr
&=1-\frac{1}{KT}\sum_{k\in \scr{K}}\sum_{t\in \scr{T}}\frac{Z}{G_k(t)+Z}
\end{align}
\end{footnotesize}

Since $G_k(t)>0$ and $Z>0$, each element in the summation is less than $1$ and overall value of the summation is less than $KT$, and this clearly leads to $1+v^T\hat{A}^{-1}u\neq0$. Note that
\begin{eqnarray*}
\hat{A}^{-1}uv^T\hat{A}^{-1}=\frac{-Z}{(KT)^2} \times \hspace{2in} \\
\begin{footnotesize} \begin{pmatrix}
\frac{1}{(G_1(1)+Z)^2}& \dots &\dots&  \frac{1}{(G_1(1)+Z)(G_K(T)+Z)}\\
\frac{1}{(G_1(1)+Z)(G_1(2)+Z)}& \frac{1}{(G_1(2)+Z)^2}&\dots&  \frac{1}{(G_1(2)+Z)(G_K(T)+Z)}\\
\vdots & \ddots\\
\frac{1}{(G_1(1)+Z)(G_K(T)+Z)}& \dots &\dots&  \frac{1}{(G_K(T)+Z)^2}
\end{pmatrix} \end{footnotesize}
\end{eqnarray*}
Thus,  \begin{align} \nonumber &A^{-1}=\begin{footnotesize}\begin{pmatrix}
\frac{1}{KT(G_1(1)+Z)}& 0&\dots& 0\\
0 & \frac{1}{KT(G_1(2)+Z)}& \dots& 0\\
\vdots & \ddots\\
0 &\dots&0 & \frac{1}{KT(G_K(T)+Z)}
\end{pmatrix} \end{footnotesize}\cr
&+\frac{Z}{(KT)^2-KT\sum_{k\in \scr{K}}\sum_{t\in \scr{T}}\frac{Z}{G_k(t)+Z}} \times \cr
&\begin{footnotesize}\begin{pmatrix}
\frac{1}{(G_1(1)+Z)^2}& \dots &\dots&  \frac{1}{(G_1(1)+Z)(G_K(T)+Z)}\\
\frac{1}{(G_1(1)+Z)(G_1(2)+Z)}& \frac{1}{(G_1(2)+Z)^2}&\dots&  \frac{1}{(G_1(2)+Z)(G_K(T)+Z)}\\
\vdots & \ddots\\
\frac{1}{(G_1(1)+Z)(G_K(T)+Z)}& \dots &\dots&  \frac{1}{(G_K(T)+Z)^2}\end{pmatrix}\end{footnotesize} \end{align}
Since $P=A^{-1}Y$, the price selection is uniquely given by \begin{equation} p_k(t)=\frac{B}{KT(G_k(t)+Z)}\left(1+\frac{\sum_{k\in \scr{K}}\sum_{t\in \scr{T}}\frac{Z}{G_k(t)+Z}}{KT-\sum_{k\in \scr{K}}\sum_{t\in \scr{T}}\frac{Z}{G_k(t)+Z}}\right)\end{equation}
which simplifies to the expression (\ref{p}).

\item Suppose that to the contrary, $p_k(t)\leq0$, this leads to
$$\frac{B}{KT(G_k(t)+Z)}\left(1+\frac{\sum_{k\in \scr{K}}\sum_{t\in \scr{T}}\frac{Z}{G_k(t)+Z}}{KT-\sum_{k\in \scr{K}}\sum_{t\in \scr{T}}\frac{Z}{G_k(t)+Z}}\right) \leq 0$$ for some $t$ and $k$.
Note that since $\frac{B}{KT(G_k(t)+Z)}$ is non-negative, there holds
$$1+\frac{\sum_{k\in \scr{K}}\sum_{t\in \scr{T}}\frac{Z}{G_k(t)+Z}}{KT-\sum_{k\in \scr{K}}\sum_{t\in \scr{T}}\frac{Z}{G_k(t)+Z}} \leq 0$$
This implies that 
$$\sum_{k\in \scr{K}}\sum_{t\in \scr{T}}\frac{Z}{G_k(t)+Z}\leq-\left(KT-\sum_{k\in \scr{K}}\sum_{t\in \scr{T}}\frac{Z}{G_k(t)+Z}\right)$$
and hence $KT\leq0$. But $K\geq1$ and $T\geq1$. Thus, this is a contradiction. Therefore, we conclude that $p_k(t)>0\,\,\,\, \forall \,t \in \scr{T}, k \in \scr{K}$.

\item Given the prices announced by other UCs, suppose that a UC $k$ announces a price of $p'_k(t)=p_k(t)+\epsilon$ at a fixed time $t$. If $\epsilon<0$, the company will decrease its price. But, the total demand from all users from UC $k$ cannot exceed the power availability (recall that the solution in 1) happens when all the available power is being sold). Since the revenue at time $t$ is the total demand multiplied by the price, UC $k$ does not achieve a higher revenue when $\epsilon<0$. Now suppose that $\epsilon>0$; this leads to
    $$U_{gen,k}(p'_k(t),{\bf{p_{-k}(t)}})-U_{gen,k}(p_k(t),{\bf{p_{-k}(t)}}) $$
    \begin{eqnarray*}
    &=&
    (p_k(t)+\epsilon) \frac{B+Z\sum_{t\in \scr{T}}p_k(t)+Z\epsilon}{KT(p_k(t)+\epsilon)} \\
    && -(p_k(t)+\epsilon)Z-p_k(t)\frac{B+Z\sum_{t\in \scr{T}}p_k(t)}{KTp_k(t)}+Zp_k(t) \\
    &=& -Z\epsilon \frac{KT-1}{KT}
    \end{eqnarray*}
    But $KT\geq1$ and $Z\geq1$. Thus, if $\epsilon>0$, the UC does not achieve a higher revenue. Therefore, for every period $t$, company $k$ does not benefit from deviating from (\ref{p}). Since this applies to every period, it applies for the entire time horizon because of the linearity of the revenue function (it is a linear combination of the demands multiplied by the prices).
\end{enumerate}
\end{proof}

In practice, due to production costs and market regulations, $p_k(t)$ cannot be outside the range of some lower and upper bounds $[p^{{\rm min}}_k(t),p^{{\rm max}}_k(t)]$  for all $t \in \scr{T}$ and $k \in \scr{K}$, as in \cite{sabita}. If $p_k(t)<p^{{\rm min}}_k(t)$, then $p_k(t)$ is set to $p^{{\rm min}}_k(t)$, and similarly for the upper-bound, if  $p_k(t)>p^{{\rm max}}_k(t)$, then we set  $p_k(t)=p^{{\rm max}}_k(t)$. The expression (\ref{p}) will still hold for the other prices because of the convexity of the problem.

\begin{remark} Using the expression (\ref{UC}), it can be verified that the revenue function for UC $k$ is

\begin{align}  &U_{{\rm gen},k}(p_k,{\bf{p_{-k}}})=\sum_{t\in \scr{T}}p_k(t)\sum_{n\in \scr{N}}d_{n,k}(t) \cr &=  \sum_{t\in \scr{T}}p_k(t)(\frac{\sum_{n\in \scr{N}}B_n+\sum_{n\in \scr{N}}\zeta_n\sum_{k\in \scr{K}}\sum_{t\in \scr{T}}p_k(t)}{KTp_k(t)} - \sum_{n\in \scr{N}}\zeta_n) \cr &=  \sum_{t\in \scr{T}}\frac{B+Z\sum_{k\in \scr{K}}\sum_{t\in \scr{T}}p_k(t)}{KT} -  Z\sum_{t\in \scr{T}}p_k(t) \cr &=  \frac{B}{K}+\frac{Z(\sum_{k\in \scr{K}}\sum_{t\in \scr{T}}p_k(t)-K\sum_{t\in \scr{T}}p_k(t))}{K} \end{align}
We can see that the revenue function of each company depends on the price selections of other companies and thus it is appropriate to use the notion of Nash equilibrium.
\end{remark}

\begin{remark}
Since $Z=\sum_{n\in \scr{N}} \zeta_n$ and $\zeta_n$ is typically 1, the value of $Z$ typically equals $N$. In this case, by (\ref{p}), we observe that for any given $G_k(t)$'s,
$$p_k(t)(G_k(t)+Z)=p_k(t)(G_k(t)+N)$$
which is  a constant for all $t \in \scr{T}$ and $k \in \scr{K}$.
Thus, the power availability is inversely proportional to the prices.
\end{remark}

\begin{remark}
Theorem 2 provides a computationally cheap expression for the prices. Since $p_k(t)$ can be directly computed using (\ref{p}), there is no need to numerically compute $A^{-1}$ or $|A|$. This enables us to deal with a large number of periods or UCs, without worrying about computational complexity.
\end{remark}
\begin{figure*}[t]
\centering
\includegraphics[width=7.5in,height=4.5in]{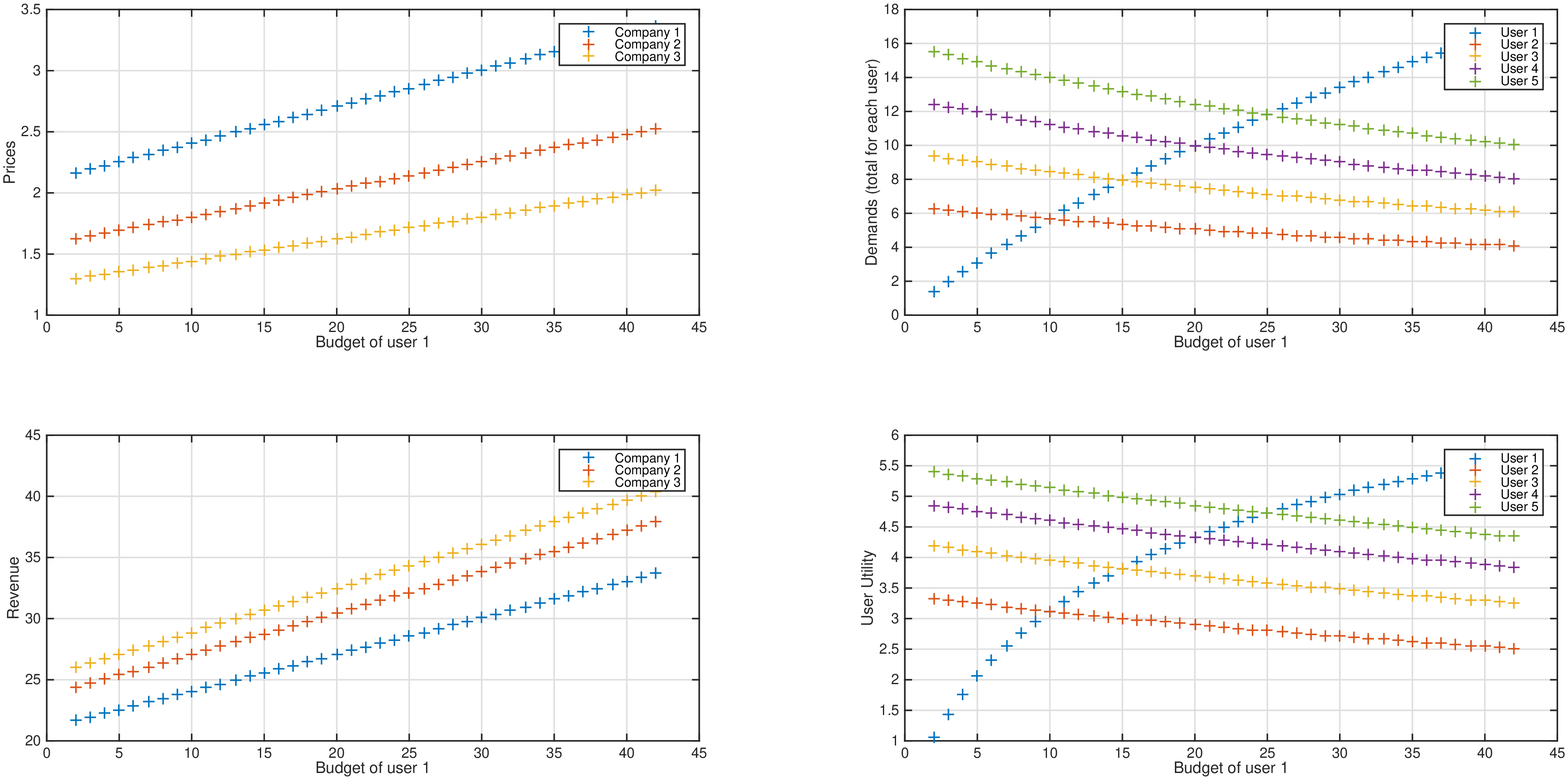}
\caption{Computations with 5 users, $T=1$,  and varying $B_1$}
\end{figure*}

\begin{figure*}[t]
\centering
\includegraphics[width=7.5in,height=4.5in]{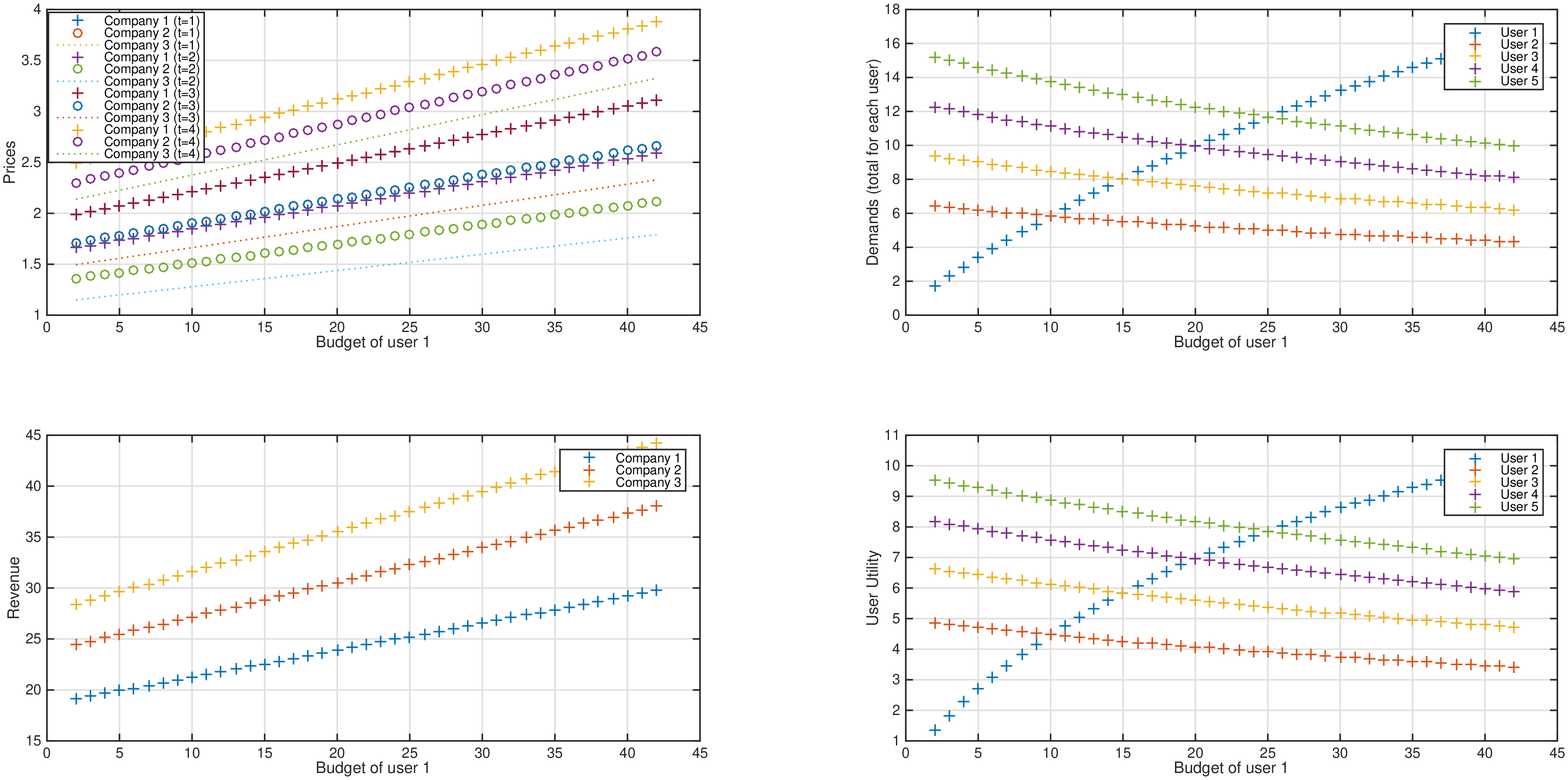}
\caption{Computations with 5 users, $T=4$,  and varying $B_1$}
\end{figure*}

\subsection {Existence and Uniqueness of the Stackelberg Equilibrium}

Companies play a noncooperative game at the market level and announce their prices to the consumers. As discussed above, consumers optimally respond to these prices and the response of user $n\in \scr{N}$ is uniquely given by (\ref{xx}). We assume that the participants of the game at the followers level have a sufficient budget to participate (derived in Theorem 1). Denote the strategy space of user (follower in the game) $n \in \scr{N}$ by $\scr{F}_n$ and the strategy space of all users by $\scr{F}= \scr{F}_1\times\dots\times\scr{F}_N$. Denote the strategy space of UC (a leader in the game) $k\in\scr{K}$ at $t\in\scr{T}$ by $\scr{L}_{k,t}:=[p^{{\rm min}}_k(t),p^{{\rm max}}_k(t)]$. Note that $p_k(t)\in\scr{L}_{k,t}$ for all $t \in \scr{T}$ and $k \in \scr{K}$. The strategy space of UC $k$ for the entire time horizon is $\scr{L}_{k}=\scr{L}_{k,1}\times\dots\times\scr{L}_{k,T}$ and the strategy space of all companies is $\scr{L}=\scr{L}_{1}\times\dots\times\scr{L}_{K}$. Before stating our main theorem, we need the following game-theoretic concepts from \cite{basar}.

The vector of prices ${\bf{p^*}} \in \scr{L}$ constitutes a {\em Nash equilibrium} for the price selection game at the UCs-level if $$U_{{\rm gen},k}(p^*_k,{\bf{p^*_{-k}}})\geq U_{{\rm gen},k}(p_k,{\bf{p_{-k}}^*}), \;\; \forall \; p_k \in \scr{L}_k$$
For given price selections  $(p_1,\dots,p_K) \in \scr{L}_1 \times \dots \times \scr{L}_K$, the optimal response from all users is
$${\bf{d^*(p)}}=\{d_1^*({\bf{p}}),d_2^*({\bf{p}}),\dots,d_N^*({\bf{p}})\}$$
where for each $n \in \scr{N}$, $d_n^*({\bf{p}})$ is the unique maximizer for $U_{{\rm user},n}(d_n,{\bf{p}})$ over $d_n \in \scr{F}_n$.
For the game considered here, the {\em Stackelberg equilibrium} is defined as $({\bf{d^*(p)}},{\bf{p^*}})$.

In general, in Stackelberg games, the response from the followers has to be unique for the equilibrium to be well defined \cite{basar}. In the game here, due to market competition, leaders aim to choose their prices in the most profitable way while taking into account what other leaders are doing. We capture the competition on the leaders' level through a Nash game. We note that in the parlance of dynamic game theory \cite{basar}, we are dealing here with open-loop information structures, with the corresponding equilibria at the utilities level being open-loop Nash equilibrium. Therefore, this is a one-shot game at which all the prices for the all periods are announced at the beginning of the game, and the followers respond to these prices by solving their local optimization problems. We have the following theorem:

\begin{theorem}
The following statements are true.
\begin{enumerate}
\item There exists a unique (open-loop) Nash equilibrium for the noncooperative game at the leaders' level.
\item There exists a unique (open-loop) Stackelberg equilibrium.
\item The maximizing demands given by (\ref{xx}) and the revenue-maximizing prices given in Theorem 2 constitute the (open-loop) Stackelberg equilibrium for the demand response management game.
\end{enumerate}
\end{theorem}

\vspace{.1in}
\begin{proof}
\begin{enumerate}
\item Note that the set $\scr{L}_{k,t}:=[p^{min}_k(t),p^{max}_k(t)] \in \mathbb{R}$ is compact (closed and bounded) $\forall \,t \in \scr{T}, k \in \scr{K}$. Thus, $\scr{L}=\scr{L}_{1}\times\dots\times\scr{L}_{K}$ is compact subset of $\mathbb{R}^{KT}$. Furthermore,
    $$\frac{\partial^2 U_{gen,k}(p_k(t),{\bf{p_{-k}}})}{\partial p^2_k(t)}=0 \,\,\,\, \forall \,t \in \scr{T}, k \in \scr{K}$$
    since the revenue function is convex and concave in $p_k(t)$.
    Hence, there exists a Nash equilibrium for the noncooperative game at the leaders level \cite{basar}. By Theorem 2, $p_k(t)$ that constitutes the best response of company $k \in \scr {K}$ to prices set by other companies is uniquely given by (\ref{p}). Thus, the Nash equilibrium is unique.
\item From 1), a unique Nash equilibrium exists at which the maximizing prices are announced to the consumers. Since the maximizing demands were uniquely given by (\ref{xx}), then there exists a unique Stackelberg equilibrium.
\item Immediately follows from Theorem 2 and parts 1)-2)\end{enumerate}
\end{proof}

\section{Numerical Results} \label{simulation}

We conducted numerical computations capturing different scenarios. In the first case, we set $T=1$ and compute the results with the same parameters values as in the single-period case \cite{sabita}. Here, for $n \in \scr{N}$, we have $\zeta_n=1$ and  $\gamma_n=1$. Additionally, $K=3$, $N=5$,  $B_1$ varies from 2 to 42, $B_2=10$,  $B_3=15$, $B_4=20$, and $B_5=25$. The power availabilities are $G_1(1)=10$, $G_2(1)=15$, and $G_3(1)=20$. Figure 1 shows that the results matches the single-period case given in \cite{sabita}. This is expected since the multi-period Stackelberg game is a generalization of the single-period one. Now, we let $T=4$, an interpretation for which can be the following: morning, afternoon, evening, and late night. The budgets of consumers are kept the same, and as before, $B_1$ still varies from 2 to 42. The total power availability for each company is also kept at the same level, but distributed across the 4 periods as follows: 25\%, 40\%, 25\%, and 10\%. Figure 2 shows that with the same total power availability for each company, and without increasing any of the consumers' budgets, the utilities for users increase significantly (almost doubles for user 1). The total demands for the users do not change, since they match the power availability. The trend for the revenues is similar to the single-period case. One key observation is that the multi-period scheme provides more incentives for consumers' participation (can be shown analytically), which is quite important and a key issue \cite{DOECOM}. Finally, we now increase the number of users and study the behaviors with a varying $T$. Here, $T$ varies from 1 to 50,  $N=50$, $K=1$, $G_1(t)=\frac{300}{T} \,\, \forall \,\, t \in \scr{T}$ (so, the power availability is equally distributed among all periods). $B_{1-10}=5$, $B_{11-20}=10$, $B_{21-30}=15$, $B_{31-40}=20$, and $B_{41-50}=25$. Figure 3 shows that the previous observation again holds, and increasing the number of the periods provides more incentives for users' participation without affecting the revenue of the company.
\begin{figure*}
\centering
\includegraphics[width=7.5in,height=2.5in]{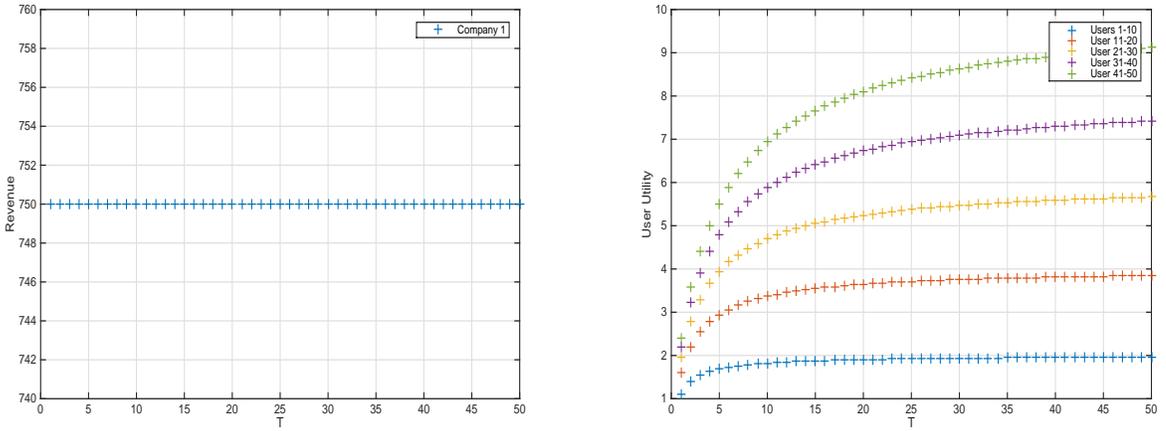}
\caption{Computations for 50 users with varying $T$}
\end{figure*}
\section{Conclusions and Future Directions} \label{conclusion}

In this paper, we have studied the multi-period DRM problem through game-theoretic methods.
In particular, we have developed a Stackelberg game to capture the interactions between UCs and energy consumers
in which UCs are the leaders and the users are the followers of the game. The UCs play a noncooperative game, which was shown to have a unique equilibrium at which companies maximize their revenues in response to price decisions made by the other companies. Then, the users optimally respond to the price selection made by UCs, by choosing the utility maximizing demand.
The overall hierarchal interaction admits a unique Stackelberg equilibrium. Furthermore, closed-form solutions have been derived, and numerical results have shown that the multi-period scheme provides more incentives for the participation of energy consumers.

There are numerous opportunities for future work. Incorporating the economic dispatch problem at the UC level is one possible direction. Another direction is to study energy scheduling and storage at the consumers level. Moreover, as this paper has focused on shiftable energy need at the consumer level, another opportunity is to study the addition of some period-specific energy need constraints.

\bibliographystyle{IEEEtran}
\bibliography{references}


\end{document}